\documentclass{article}
\begin{document}

\newtheorem{lemma}{Lemma}
\newtheorem{theorem}{Theorem}
\newtheorem{definition}{Definition}
\newtheorem{corollary}{Corollary}
\newtheorem{proposition}{Proposition}

\def\R{{\bf R}}
\def\C{{\bf C}}
\def\Z{{\bf Z}}
\def\Q{{\bf Q}}
\def\A{{\cal A}}
\def\M{{\cal M}}
\def\vol{\mbox{vol}}

\begin{center}
{\Large \bf Integrable geodesic flow with positive topological
entropy}
\end{center}
\begin{center} {Alexey V. BOLSINOV
\footnote{
    Department of Mathematics and Mechanics,
    Moscow State University, 119899 Moscow, Russia, e-mail:
    bols@difgeo.math.msu.su}
and Iskander A. TAIMANOV
\footnote{
    Institute of Mathematics,
    630090 Novosibirsk, Russia,
    e-mail: taimanov@math.nsc.ru}}
\end{center}

\bigskip

\section{Introduction and main results}

The main result of this paper is the following theorem.

\begin{theorem}
There is a real-analytic Riemannian manifold $M_A$ diffeomorphic to
the quotient of $T^2 \times
{\R}^1$ with respect to the free ${\Z}$-action generated
by the map
$$
(X,z) \to (AX,z+1),
$$
where $X = (x,y) \in T^2 = {\R}^2/{\Z}^2$,
$z \in {\R}$, and $A$ is the Anosov automorphism of
the $2$-torus $T^2$ defined by the matrix
$$
A = \left( \begin{array}{cc} 2 & 1 \\ 1
& 1 \end{array} \right),
\eqno{(1)}
$$
such that

i) the geodesic flow on $M_A$ is (Liouville) integrable
by $C^\infty$ first integrals;

ii) the geodesic flow on $M_A$ is not (Liouville) integrable
by real-analytic first integrals;

iii) the topological entropy of the geodesic flow $F_t$ is
positive;

iv) the fundamental group $\pi_1(M_A)$ of the manifold $M_A$
has an exponential growth;

v) the unit covector bundle $S M_A$ contains a submanifold
$N$ such that $N$ is diffeomorphic to the $2$-torus $T^2$
and the restriction of $F_1$ onto $N$ is the
Anosov automorphism given by matrix (1).
\end{theorem}

To explain the statement in detail we recall main definitions
and results on topological obstructions to integrability of
geodesic flows.

Let $g_{jk}$ be a Riemannian metric on an
$n$-dimensional manifold $M^n$.
It defines the geodesic flow on the tangent bundle $T M^n$
which is a Lagrangian system with the Lagrange function
$$
L(x,\dot{x}) = \frac{1}{2} g_{jk} \dot{x}^j \dot{x}^k.
$$
The Legendre transform $T M^n \to T^\ast M^n$
$$
\dot{x}\in T_x M^n \to p \in T_x^\ast M^n: \ p_j =
g_{jk}\dot{x}^k
$$
maps this Lagrangian system into a
Hamiltonian system on $T^\ast M^n$
with a symplectic form
$$
\omega = \sum_{j=1}^n d x^j \wedge d p_j
$$
and the Hamilton function
$$
H(x,p) = \frac{1}{2}g^{jk}(x) p_j p_k.
$$
This Hamiltonian system is also called the geodesic flow of the
metric.

The symplectic form defines the Poisson brackets on the space of
smooth functions on $T^\ast M^n$ by the formula
$$
\{ f, g\} = h^{jk} \frac{\partial f}{\partial y^j}
\frac{\partial g}{\partial y^k},
$$
where $\omega = h_{jk} d y^j \wedge d y^k$ locally.

It is said that a Hamiltonian system on a $2n$-dimensional
symplectic manifold is (Liouville) integrable if there are
$n$ first integrals $I_1,\dots,I_n$ of this system such that

1) these integrals are in involution: $\{I_j, I_k\} = 0$ for any
$j,k$, $1 \leq j, k \leq n$;

2) these integrals are functionally independent almost
everywhere, i.e., on a dense open subset.

Since restrictions of the geodesic flow onto different
non-zero level surfaces of its Hamilton function $H=I_n$ are
smoothly trajectory equivalent, we may replace this notion of
integrability by the weaker condition that there are $n-1$
additional first integrals $I_1,\dots,I_{(n-1)}$ which are in
involution and functionally independent almost everywhere on the
unit covector bundle $S M^n = \{H(x,p) = 1\} \subset T^\ast M^n$
\cite{T1}.

If $M^n$ is real-analytic together with the metric and all first
integrals $I_1,\dots,I_n$, then it is said that the geodesic
flow is analytically integrable.

Kozlov established the first
topological obstruction to analytic integrability of geodesic
flows proving that the geodesic flow of a real-analytic metric
on a two-dimensional closed oriented manifold $M^2$ of genus
$g > 1$ does not admit an additional analytic first integral
\cite{K1} (see also \cite{K2} for general setup of
nonintegrability problem).

For higher-dimensional manifolds, obstructions to integrability
were found in \cite{T1,T2} where it was proved that
analytic integrability of the geodesic flow on a manifold $M^n$
implies that

1) the fundamental group $\pi_1(M^n)$ of $M^n$ is almost
commutative, i.e., contains a commutative subgroup of finite
index;

2) if the first Betti number $b_1(M^n)$ equals $k$, then
the real cohomology ring $H^\ast(M^n;{\R})$ contains a
subring isomorphic to the real cohomology ring of a
$k$-dimensional torus and, in particular, $b_1(M^n) \leq \dim
M^n = n$.

In these results the analyticity condition may be
replaced by stronger geometric condition called geometric
simplicity and reflecting some tameness properties of the
singular set where the first integrals are functionally
dependent. For instance, one may only assume that $S M^n$ is a
disjoint union of a closed invariant set $\Gamma$ which is
nowhere dense and of finitely many open toroidal domains
foliated by invariant tori.

Later Paternain proposed another approach to finding topological
obstructions to integrability based on a vanishing of the
topological entropy of the geodesic flow on $S M^n$. If
this quantity vanishes, then $\pi_1(M^n)$ has a subexponential
growth \cite{Din} and, if in addition $M^n$ is a $C^\infty$
simply-connected manifold, then $Y$ is rationally-elliptic
(this follows from results of Gromov and Yomdin) \cite{P1,P2}.
Integrability implies vanishing of the topological entropy under
some additional conditions which were established in
\cite{P1,P2,T3} and restrict not only the singular
set but also the behaviour of the flow on this set.

Recently Butler has found new examples of $C^\infty$ integrable
geodesic flows of homogeneous metrics on nilmanifolds
\cite{Butler}. The
simplest of them is a $3$-manifold $M_B$ obtained from a
product $T^2 \times [0,1]$ by
identifying the components of the boundary by a homeomorphism
$$
(X,0) \to (BX,1),
$$
where
$$
B = \left( \begin{array}{cc}
1 & 1 \\ 0 & 1 \end{array} \right)
$$
and $X \in T^2$.
The fundamental group of the resulting manifold $M_B$ is not
almost commutative, $b_1(M_B) = 2$, and
$H^\ast(M_B;{\R})$ does not contain a subring isomorphic
to $H^\ast(T^2;{\R})$. However, $b_1(M_B) < \dim
M_B$.  This shows that some of the results of \cite{T1,T2} are
not generalized for the $C^\infty$ case.  Note that the
topological entropy vanishes for Butler's examples.

The present paper is based on an observation that Butler's
construction is generalized for constructing $C^\infty$
integrable geodesic flows on all $T^n$-bundles over $S^1$.
In the case when the gluing automorphism $C:T^n \to T^n$ is
hyperbolic we obtain remarkable Hamiltonian systems on a
cotangent bundle to $M_C$:  they are $C^\infty$ integrable but
have positive topological entropy. This, in particular, shows
that treating positivity of topological entropy as a criterion
for chaos which is used sometimes is not correct.

We confine only to one example of such a flow which we study in
detail.

\section{The metric on $M_A$ and its geodesic flow}

Let
$$
A:{\Z}^2 \to {\Z}^2
$$
be an automorphism determined by matrix (1).
It determines the following action on $T^2$:
$$
(x,y) \, \mbox{mod}\, {\Z}^2
\to (2x +y, x + y)\, \mbox{mod}\, {\Z}^2.
$$

We construct $M_A$ as follows. Take a product $T^2 \times
[0,1]$ and identify the components of its boundary
using the automorphism $A$:
$$
(X,0) \sim (AX,1),
$$
where $X = (x,y) \in T^2$. We denote the resulted manifold by
$M_A$.  Near every point $p \in M_A$ we have local coordinates
$x, y$ and $t$, where $z$ is a linear coordinate on
$S^1 = {\R}/{\Z}$.

Take the following metric on $M_A$:
$$
d s^2 = d z^2 + g_{11}(z) d x^2 + 2g_{12}(z) d x d y +
g_{22}(z) d y^2
$$
where
$$
G(t) =
\left( \begin{array}{cc} g_{11} & g_{12} \\
g_{21} (= g_{12}) & g_{22} \end{array} \right) =
\exp(-z G_0^\top) \exp(-z G_0)
\eqno{(2)}
$$
and $\exp G_0 = A$. We set
$$
g_{33}=1, g_{13}=g_{23} = 1.
$$
Indeed, this formula defines a metric on an infinite cylinder
${\cal C} = T^2 \times {\R}$ which is invariant with respect
to the ${\Z}$-action generated by
$$
(x,y,z) \to (2x+y,x+y,z+1),
\eqno{(3)}
$$
and, therefore, it descends to a metric on the quotient space
$M_A = {\cal C}/{\Z}$.

\begin{proposition}
The geodesic flow of metric (2) on the infinite
cylinder ${\cal C}$ admits three first integrals which are
functionally independent almost everywhere.
\end{proposition}

{\sl Proof of Proposition.}
The Hamiltonian function
$$
F_3 = H = \frac{1}{2}\left(
p_z^2 + g^{11}(z)p_x^2 + 2g^{12}(z)p_x p_y + g^{22}(z) p_y^2
\right)
$$
of this flow is, by the construction, a first integral.
Since $H$ does not depend on $x$ and $y$ the quasimomenta
$F_1 = p_x$ and $F_2 = p_y$ are also first integrals.
It is clear that the set of first integrals $I_1,I_2$, and
$I_3$ is functionally independent almost everywhere.
This proves the proposition.

Since action (3) preserves the symplectic form $\omega$, it
induces the following action on $T^\ast {\cal C}$:
$$
\left(
\begin{array}{c}
p_x \\ p_y
\end{array}
\right)
\to
\left(
\begin{array}{cc}
1 & -1 \\ -1 & 2
\end{array}
\right)
\left(
\begin{array}{c}
p_x \\ p_y
\end{array}
\right),
\ \
p_z \to p_z.
$$
This descends to a linear action on $T^\ast M_A$
which preserves fibers and takes the form
$$
\left( p_x - \frac{1+\sqrt{5}}{2}p_y \right)
\to \lambda
\left( p_x - \frac{1+\sqrt{5}}{2}p_y \right),
$$
$$
\left( p_x - \frac{1-\sqrt{5}}{2}p_y \right)
\to
\lambda^{-1} \left( p_x - \frac{1-\sqrt{5}}{2}p_y \right),
\eqno{(4)}
$$
$$
p_z \to p_z,
\ \ \lambda =
\frac{3+\sqrt{5}}{2}.
$$
It is evident that the indefinite quadratic form
$$
I_1 = \left( p_x -
\frac{1+\sqrt{5}}{2}p_y \right) \left( p_x -
\frac{1-\sqrt{5}}{2}p_y \right) = p_x^2 -p_x p_y - p_y^2
$$
and the positively definite quadratic form
$$
I_3 = H = \frac{1}{2}\left(
p_z^2 + g^{11}(z)p_x^2 + 2g^{12}(z)p_x p_y + g^{22}(z) p_y^2
\right)
$$
are invariants of this action.
To construct the third invariant we notice that
$$
\frac{\log \big|p_x - \frac{1+\sqrt{5}}{2}p_y \big|}
{\log \lambda}
$$
is not invariant but the action adds $1$ to this
quantity when it is correctly defined. Therefore, the following
function
$$
I_2 = f(I_1)\cdot \sin \left( \frac{\log \big|p_x -
\frac{1+\sqrt{5}}{2}p_y \big|}{\log \lambda} \right),
$$
where
$$
f(u) = \exp\left(-\frac{1}{u^2}\right),
$$
is everywhere defined and invariant with respect to action (4).

\begin{proposition}
The functions $I_1, I_2$, and $I_3$ are $C^\infty$ first
integrals of the geodesic flow on $M_A$ which are
functionally independent almost everywhere. Therefore,
the geodesic flow on $M_A$ is (Liouville) integrable by
$C^\infty$ functions.
\end{proposition}

{\sl Proof of Proposition.}
The functions $I_1, I_2$ , and $I_3$ on $T^\ast {\cal C}$ are
invariants of action (4) and, therefore, descend to functions on
$T^\ast M_A$. We may consider $I_1$ and $I_2$ as replacing $F_1$
and $F_2$: they are pairwise involutive and independent on
spatial variables $x,y,z$. Moreover they do not depend on $p_z$
and, therefore, they  are in involution with $I_3=H$ which, in
particular, means that they are first integrals of the geodesic
flow. It remains to notice that, by their construction, they
are $C^\infty$. This finishes the proof of Proposition.

\begin{proposition}
Let $N$ be a subset of the unit covector bundle $S M_A$ formed
by the points with
$$
z = 0, \ \ p_x = p_y = 0, \ \ p_z = 1.
$$
Then it is diffeomorphic to $T^2$ and the translation
$$
F_t: T^\ast M_A \to T^\ast M_A
$$
along the trajectories of the geodesic flow for $t=1$ maps $N$
into itself and this map is the Anosov
automorphism given by matrix (1).
\end{proposition}

{\sl Proof of Proposition.}
The geodesic flow on $M_A$ is covered by the geodesic flow on
${\cal C}$ for which $p_x$ and $p_y$ are first integrals.
Therefore, on ${\cal C}$ the translation of the preimage of
$N$ under projection is as follows:
$$
F_t(x,y,z,p_x,p_y,p_z) = F_t(x,y,0,0,0,1) =
(x,y,t,0,0,1).
$$
Recalling the  construction of $M_A$ proves the proposition.

Note that Propositions 2 and 3 prove statements i) and v)
of  Theorem 1, respectively.

\section{The fundamental group $\pi_1(M_A)$ and the topological
entropy of the geodesic flow on $M_A$}

The manifold $M_A$ is covered by ${\R}^3$ on which acts
$\pi_1(M_A)$. This group is generated by
$$
a:(x,y,z) \to (x+1,y,z), \ \
b:(x,y,z) \to (x,y+1,z),
$$
$$
c:(x,y,z) \to (2x+y,x+y,z+1).
$$
The relations between these generators are
$$
[a,b] = 1, \ \ [c,a] = ab, \ \ [c,b] = a.
$$

\begin{proposition} [see, for instance, \cite{GK}]
$\pi_1(M_A)$ has an exponential growth.
\end{proposition}

This follows from the hyperbolicity of $A$ or may be proved
directly: the words $c a^{\varepsilon_1} c a^{\varepsilon_2}
\dots c a^{\varepsilon_k}$ are different for
$\varepsilon_j = 0,1$ and, therefore, $\gamma(2k) \geq 2^k$,
where $\gamma$ is the
growth function of $\pi_1(M_A)$ with respect to generators
$a,b$, and $c$.

\begin{corollary}
The geodesic flow on $M_A$ is not (Liouville) integrable by
real-analytic first integrals.
\end{corollary}

It follows from the results of \cite{T1} (also exposed
in Section 1) that if this flow is analytically integrable, then
$\pi_1(M_A)$ is almost commutative and, therefore, has
a polynomial growth. This contradiction establishes the
corollary.

\begin{corollary}
The topological entropy of the geodesic flow on $M_A$ is
positive.
\end{corollary}

Indeed, it was proved by Dinaburg, that if the fundamental group
of a manifold $M^n$ has an exponential growth, then the
topological entropy of the geodesic flow of any Riemannian
metric on $M^n$ is positive \cite{Din}.

The latter corollary also follows from Proposition 3: it is
known that the topological entropy equals the supremum of
the measure entropies taken over all ergodic invariant Borel
measures. Hence, we may take a singular measure concentrated on
$N \subset S M_A$ which has the form
$$
d \mu = d x \wedge d y.
$$
It is well known that the measure entropy of the
Anosov automorphism $A:N \to N$ is positive (this follows,
for instance, from nonvanishing of the Lyapunov exponents for
any point of $N$).

Note that Proposition 4 and Corollaries 1 and 2 establish
statements iv), ii), and iii) of Theorem 1, respectively.

\bigskip

The authors thank L. Butler for sending them his preprint
\cite{Butler} and I. K. Babenko for helpful discussions.

The authors were partially supported by Russian Foundation for
Basic Researches (grants 96-15-96868 and 98-01-00240 (A. V.
B.), and 96-15-96877 and 98-01-00749 (I.A.T.)) and INTAS (grant
96-0070 (I.A.T.)).

\end{document}